\documentclass[runningheads]{llncs}
\usepackage{graphicx}
\usepackage{tikz}
\usepackage{cite}
\usepackage{amsmath,amssymb,amsfonts}
\usepackage{algorithmic}
\usepackage{graphicx}
\usepackage{textcomp}
\usepackage{xcolor}

\usepackage{amsmath}

\usepackage{hpmacros04}

\usepackage{xcolor}
\colorlet{darkmagenta}{magenta!85!black}
\colorlet{darkgreen}{green!55!black}
\colorlet{darkblue}{blue!85!black}
\colorlet{darkred}{red!80!black}

\newcommand{\al}{\ensuremath{\alpha}}
\newcommand{\be}{\ensuremath{\beta}}

\newcommand{\tc}[1]{\ensuremath{\mathbf{#1}}} 
\newcommand{\tct}{\tc t} 
\newcommand{\tcf}{\tc f} 
\newcommand{\tcb}{\tc b} 
\newcommand{\tcn}{\tc n} 

\newcommand{\At}{\ensuremath{\mathit{At}}}
\newcommand{\Var}{\ensuremath{\mathit{Var}}}
\newcommand{\Const}{\ensuremath{\mathit{Const}}}
\newcommand{\Pred}{\ensuremath{\mathit{Pred}}}

\newcommand{\tle}{\ensuremath{\mathbin{\le_\mathrm t}}}
\newcommand{\ile}{\ensuremath{\mathbin{\le_\mathrm i}}}

\newcommand{\model}[1]{\ensuremath{\boldsymbol{#1}}}
\newcommand{\classmod}[1]{\ensuremath{\boldsymbol{#1}}}

\newcommand{\exP}{\ensuremath{\mathrm{E!}}} 
\newcommand{\inD}{\ensuremath{D_1}} 
\newcommand{\outD}{\ensuremath{D_0}} 
\newcommand{\inAll}{\all}
\newcommand{\inExi}{\exi}
\newcommand{\outAll}{\mathrm{\Pi}}
\newcommand{\outExi}{\mathrm{\Sigma}}

\newcommand{\SemVal}[3]{\ensuremath{\semvalue{#1}_{\model{#2},{#3}}}}
\newcommand{\SemValPlus}[3]
	{\ensuremath{\semvalue{#1}\!^+_{\model{#2},{#3}}}}
\newcommand{\SemValMinus}[3]
	{\ensuremath{\semvalue{#1}\!^-_{\model{#2},{#3}}}}

\newcommand{\logic}[1]{\ensuremath{\mathrm{#1}}}
\newcommand{\BD}{\logic{BD}}
\newcommand{\BDD}{\logic{BD\Delta}}
\newcommand{\LukD}{\logic{\Luk_\baaz}}

\newcommand{\LBDD}{\logic{\Luk\BDD}}

\renewcommand{\baaz}{\ensuremath{\triangle}}
\renewcommand{\disj}{\ensuremath{\oplus}}
\newcommand{\bddelta}{\ensuremath{\boldsymbol\delta}} 
\newcommand{\bdnega}{\ensuremath{{\sim}}}
\newcommand{\bdcons}{\ensuremath{{\circ}}}

\newcommand{\bdnegad}{\ensuremath{\neg}}

\DeclareMathOperator{\sgn}{sgn}

\newcommand{\und}{\ensuremath{{*}}}

\begin{document}
\title{Free Quantification in Four-Valued and Fuzzy Bilattice-Valued Logics}
%
%
\author{Libor B{\v{e}}hounek\inst{1}\orcidID{0000-0001-8570-9657} \and
Martina Da\v{n}kov{\'a}\inst{1}\orcidID{0000-0001-5806-7898} \and
Anton{\'{i}}n Dvo{\v{r}}{\'{a}}k\inst{1}\orcidID{0000-0002-0611-491X} }
\authorrunning{L. B{\v{e}}hounek et al.}
%
\institute{University of Ostrava, IRAFM, 30.~dubna~22, 701\,03 Ostrava, Czech Republic \\
\email{\{libor.behounek|martina.dankova|antonin.dvorak\}@osu.cz}}
\maketitle              
\begin{abstract}
We introduce a variant of free logic (i.e., a logic admitting terms with nonexistent referents) that accommodates truth-value gluts as well as gaps. Employing a suitable expansion of the Belnap--Dunn four-valued logic, we specify a dual-domain semantics for free logic, in which propositions containing non-denoting terms can be true, false, neither true nor false, or both true and false. In each model, the dual domain semantics separates existing and non-existing objects into two subdomains, making it possible to quantify either over all objects or existing objects only. We also outline a fuzzy variant of the dual-domain semantics, accommodating non-denoting terms in fuzzy contexts that can be partially indeterminate or inconsistent.

\end{abstract}
\allowdisplaybreaks

\section{Introduction}
\label{s:intro}

The four-valued bilattice logic \BD\ (Belnap--Dunn logic, also known as \logic{FDE}) has been designed to deal with underdetermined (i.e., neither true nor false) as well as overdetermined (i.e., both true and false) modes of truth and falsity besides the usual truth values \emph{true} and \emph{false} \cite{Belnap1977Useful,Dunn:FDE1977,Omori:40FDE}. Its generalizations to $[0,1]^2$-valued logics have been used for reasoning under indeterminacy and uncertainty \cite{Klein-Majer-RafieeRad:ProbabilitiesGapsGluts,Bilkova-et-al:QualitativeTwoLayered,Bilkova-et-al:TwoLayeredParaconsistentProbabilities}. Here we employ a suitable expansion of the first-order logic~\BD\ (in particular, the logic~\BDD\ of~\cite{SanoOmori:BDDelta2014}) and its $[0,1]^2$-valued generalization over {\L}ukasiewicz fuzzy logic~\LukD\ as the background machinery to develop logics of the so-called free quantification in the presence of (either crisp or graded) indeterminacy and overdeterminacy of truth and falsity.

\emph{Free logic,} or `logic free of existential assumptions'~\cite{Bencivenga:FreeLogics-HPL2}, is a field of logic that studies quantification over terms that may be non-denoting. Such terms naturally occur in various situations, e.g., in logical analysis of natural language (definite descriptions), fictional discourse, and predicate modal logics (incl.\ epistemic and temporal). When terms over which quantification is carried include non-denoting ones, some standard quantification laws, such as universal instantiation (also known as specification) and its existential dual, are no longer valid. Naturally, the invalidity of these laws should be reflected in the semantics of free logics.

Free logics are developed in several variants (see, e.g.,~\cite{Nolt:FreeLogic-SEP}).
One of intuitively appealing variants, which we choose for our endeavor, is the so-called \emph{positive} free logic with \emph{dual-domain semantics.} Positive free logics allow atomic formulas with non-denoting terms to be true, false, or indeterminate (truth-valueless). In the dual-domain semantics, a model is equipped with two domains: the \emph{inner domain}~$D_1$ (which is the extension of the \emph{existence predicate,} traditionally denoted by~\exP) that collects existing objects, and the outer domain $D_0 \supseteq D_1$ containing all objects, both existing and non-existing. Then, the universal and existential quantifiers over existing objects can be defined in such a way that the range of variables is limited to the inner domain\cite{Carnieli_FreeLFL1}; then, however, the objects from $D_0 \setminus D_1$ that are not associated with a logical constant are unreachable by formulas of the given logic. A more flexible possibility is to use the so-called \emph{inner} and \emph{outer quantifiers.} The outer quantifiers range over~$D_0$ and are simply the standard quantifiers of the given logic. The inner quantifiers ranging over $D_1$ can be defined by relativization of the outer ones.

Since positive free logic allows formulas containing non-denoting terms to be truth-valueless, it is expedient to evaluate formulas in a logic which accommodates \emph{truth-value gaps,} such as three-valued strong Kleene logic~\logic{K_3}  \cite{Kleene:NotationOrdinal} that has the truth values \emph{true, false,} and \emph{neither.} An option less often considered is using a logic that also accommodates truth-value \emph{gluts,} or the truth values for propositions that are \emph{both} true and false. The need for truth-value gluts in positive free logic follows from the fact that nonexistent objects can have contradictory properties: e.g., a \emph{square circle} is both round (being a circle) and not round (being square). A recent system of free logic by Carnielli  and Antunes~\cite{Carnieli_FreeLFL1} does consider truth-value gluts (though not truth-value gaps), using the three-valued logic \logic{LFI1} with the truth values \emph{true, false,} and \emph{both.}
In this paper we propose a four-valued free logic that accommodates truth-value gaps as well as gluts, using the four values \emph{true, false, neither,} and \emph{both.} To this end, we employ the first-order logic~\BDD, which was introduced in~\cite{SanoOmori:BDDelta2014}. We recall the logic~\BDD\ in Section~\ref{s:BDD} and introduce the (positive dual-domain) free logic over~\BDD\ in Section~\ref{s:free-BDD}.

A fuzzy variant of a positive free logic with dual-domain semantics has been outlined in~\cite{Behounek-Dvorak:NondenotingTermsFuzzy-EUSFLAT}. The need for fuzzified free logic is rather natural, as non-denoting terms and terms denoting nonexistent objects can be encountered in fuzzy contexts just like in crisp contexts (e.g., when a fuzzy property is predicated of a nonexistent individual, as in ``Sherlock Holmes is clever''). In~\cite{Behounek-Dvorak:NondenotingTermsFuzzy-EUSFLAT}, free fuzzy logic is built over partial fuzzy logic~\cite{Behounek-Novak:TowardsFuzzyPartial-ISMVL}, where the additional truth value \und\ represents a truth-value gap. This contribution extends the approach outlined above in such a way that it incorporates also truth-value gluts, and moreover admits the underdeterminacy or contradictoriness of truth values to be graded, i.e., just partially indeterminate (like in interval-valued fuzzy logic) or partially inconsistent (when the degrees of truth and falsity sum up to more than the value~$1$). We sketch the fuzzification of the four-valued logic~\BDD\ via {\L}ukasiewicz fuzzy logic~\LukD\ and present the dual-domain semantics over the resulting logic~\LBDD\ in Section~\ref{s:fuzzy}. 
Some topics for future work are mentioned in Section~\ref{s:conclusion}.

\newpage

\section{Four-Valued Bilattice Logic}
\label{s:BDD}

The paradigmatic four-valued logic for reasoning about propositions that can be true, false, both, or neither, is the Belnap--Dunn logic~\BD~\cite{Belnap1977Useful,Dunn:FDE1977,Omori:40FDE}. However, since the logic \BD\ is expressively rather poor, it is expedient to use a suitable extension thereof. For our purposes, the logic~\BDD\ of~\cite{SanoOmori:BDDelta2014} is a suitable choice, as it contains (as definable connectives) the normality indicator as well as an implication needed for relativized quantifiers in Section~\ref{s:free-BDD}. Furthermore, the logic~\BDD\ is defined as a first-order four-valued logic with equality and has a sound and complete natural-deduction axiomatic system~\cite{SanoOmori:BDDelta2014}.

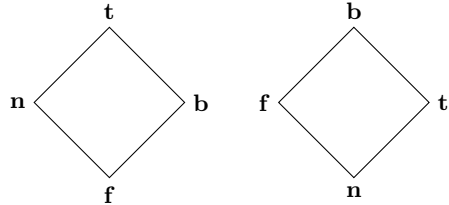
\begin{figure}
    \centering
    \begin{tikzpicture}[scale=0.5]
    \coordinate (t) at (0,2);
    \coordinate (n) at (2,0);
    \coordinate (f) at (0,-2);
    \coordinate (b) at (-2,0);

    \draw (n) -- (t) -- (b) -- (f) -- cycle;

  \node at (t) [above] {\tct};
  \node at (b) [left] {\tcn};
  \node at (f) [below] {\tcf};
  \node at (n) [right] {\tcb};

    \coordinate (t2) at (6.5,2);
    \coordinate (b2) at (8.5,0);
    \coordinate (f2) at (6.5,-2);
    \coordinate (n2) at (4.5,0);

    \draw (n2) -- (t2) -- (b2) -- (f2) -- cycle;

  \node at (t2) [above] {\tcb};
  \node at (b2) [right] {\tct};
  \node at (f2) [below] {\tcn};
  \node at (n2) [left] {\tcf};

    \end{tikzpicture}
    \caption{The Hasse diagrams of the truth order (left) and the information order (right) on the four truth values.}
    \label{fig:Orders}
\end{figure}

The set of truth values of the logic~\BDD\ is $\fclass{0,1}^2$, i.e., the set of pairs of the classical truth values $0$ and $1$. The first component~\al\ of a truth value
$\tuple{\al,\be}\in\fclass{0,1}^2$ 
indicates whether the proposition is true or not and the second component~\be\ indicates whether it is false or not; the truth and falsity of a proposition are evaluated independently, so a proposition can be both true and false as well as neither true nor false.
We can thus define the following four truth values of~\BDD:
\begin{equation*}
    \begin{array}{ccl}
     \tct=\tuple{1,0} & \quad &\text{true (only)}\\  
     \tcf=\tuple{0,1} & \quad &\text{false (only)}\\  
     \tcn=\tuple{0,0} & \quad &\text{neither true nor false}\\  
     \tcb=\tuple{1,1} & \quad &\text{both true and false}\\  
    \end{array}
\end{equation*}
It is customary to define two lattice orders on the four truth values (see Figure~\ref{fig:Orders}):
\begin{itemize}
\item
    The \emph{truth order}~\tle, where \tcf\ is the least, \tct\ the largest, and \tcn\ and \tcb\ are intermediate and mutually incomparable; and
\item 
    The \emph{information order}~\ile, where \tcn\ is the least, \tcb\ the largest, and \tcf\ and \tct\ are intermediate and mutually incomparable.
\end{itemize}
The \emph{designated} truth values (i.e., those considered ``true'' in the definition of entailment) are those which are ``at least true'' in the information order, i.e.,~\tct\ and~\tcb.
Furthermore, we call the truth values~\tct\ and~\tcf\ \emph{normal.}

Accordingly, the \emph{valuation} of propositional atoms in~\BDD\ is a function 
$v\colon\At\to\{0,1\}^2$. It can be decomposed into a pair of functions,
$v^+,v^-\colon\At\to\{0,1\}$, called the \emph{positive} and \emph{negative valuation,}
where
$v(p)=\tuple{v^+(p),v^-(p)}$
for any $p\in\At$.

The propositional language of~\BDD\ consists of the symbols for propositional connectives $\{\bdnega,\allowbreak\minc,\maxd,\bddelta\}$ with the following meanings:
\begin{center}
\begin{tabular}{c@{\qquad}l}
     \minc & (truth-order lattice) conjunction\\  
     \maxd &  (truth-order lattice) disjunction\\  
     \bdnega & negation\\  
     \bddelta & indicator of designated values\\
\end{tabular}
\end{center}
Note that we use the symbol~\bddelta\ for the indicator of designated values (whereas~\cite{SanoOmori:BDDelta2014} used $\mathrm\Delta$), to avoid confusion with the connective~\baaz\ of fuzzy logic.

The semantics of these connectives in~\BDD\ is given by the following Tarski conditions (in the form suitable for easy fuzzification):
\begin{align}
\label{eq:BDD-conn-Tarski-first}
  v^+(\f\minc\p)&=\min(v^+(\f),v^+(\p))\\   
  v^-(\f\minc\p)&=\max(v^-(\f),v^-(\p))
\\[1ex]
  v^+(\f\maxd\p)&=\max(v^+(\f),v^+(\p))\\
  v^-(\f\maxd\p)&=\min(v^-(\f),v^-(\p))
\\[1ex]
  v^+(\bdnega\f)&=v^-(\f)\\
  v^-(\bdnega\f)&=v^+(\f)
\\[1ex]
  v^+(\bddelta\f)&=v^+(\f)\\
\label{eq:BDD-conn-Tarski-last}
  v^-(\bddelta\f)&=1-v^+(\f)
  \end{align}
Thus, e.g., conjunction is (at least) true if both conjuncts are (at least) true and is (at least) false if either conjunct is (at least) false, and similarly for the other connectives (where ``at least'' means the information order~\ile).
Furthermore we define the following derived connectives of~\BDD\ (cf.~\cite{Avron:LogBilattices1996,SanoOmori:BDDelta2014,Omori:40FDE}):
\begin{align}
\label{eq:def-bdnegad-BDD}
    \bdnegad\f &\equiv
        \bdnega\bddelta\f
        &&\text{bivalent negation}
\\\label{eq:def-impl-BDD}
    \f\impl\p &\equiv
        \bdnega\f\maxd\p
        &&\text{(material) implication}
\\\label{eq:def-bdcons-BDD}
    \bdcons\f &\equiv
        (\bddelta\f\maxd\bddelta\bdnega\f)\minc
        (\bdnega\bddelta\f\maxd\bdnega\bddelta\bdnega\f)
        &&\text{normality indicator}
\end{align}

The definitions can be summarized by the following truth tables:

\begin{equation*}
\setlength{\arraycolsep}{0.7ex}
\begin{array}{r|cccc}
    \f\minc\p    & \tct & \tcb & \tcn  & \tcf
\\\hline
    \tct        & \tct & \tcb & \tcn  & \tcf
\\  \tcb        & \tcb & \tcb  & \tcf  & \tcf
\\  \tcn        & \tcn & \tcf & \tcn  & \tcf
\\  \tcf        & \tcf & \tcf & \tcf  & \tcf   
\end{array}
\qquad
\begin{array}{r|cccc}
    \f\maxd\p    & \tct & \tcb & \tcn  & \tcf
\\\hline
    \tct        & \tct & \tct & \tct  & \tct
\\  \tcb        & \tct & \tcb  & \tct  & \tcb
\\  \tcn        & \tct & \tct & \tcn  & \tcn
\\  \tcf        & \tct & \tcb & \tcn  & \tcf   
\end{array}
\qquad
\begin{array}{r|c}
    \f    & \bdnega\f
\\\hline
    \tct        & \tcf
\\  \tcb        & \tcb
\\  \tcn        & \tcn
\\  \tcf        & \tct
\end{array}
\qquad
\begin{array}{r|c}
    \f    & \bddelta\f
\\\hline
    \tct        & \tct 
\\  \tcb        & \tct 
\\  \tcn        & \tcf 
\\  \tcf        & \tcf
\end{array}
\end{equation*}

\begin{equation*}
\setlength{\arraycolsep}{0.7ex}
\begin{array}{r|cccc}
    \f\impl\p    & \tct & \tcb & \tcn  & \tcf
\\\hline
    \tct        & \tct & \tcb & \tcn  & \tcf
\\  \tcb        & \tct & \tcb  & \tct  & \tcb
\\  \tcn        & \tct & \tct & \tcn  & \tcn
\\  \tcf        & \tct & \tct & \tct  & \tct   
\end{array}
\qquad
\begin{array}{c|c}
    \f    & \bdnegad\f
\\\hline
    \tct        & \tcf
\\  \tcb        & \tcf
\\  \tcn        & \tct
\\  \tcf        & \tct
\end{array}
\qquad
\begin{array}{c|c}
    \f    & \bdcons\f
\\\hline
    \tct        & \tct
\\  \tcb        & \tcf
\\  \tcn        & \tcf
\\  \tcf        & \tct
\end{array}
\end{equation*}

The first-order language of \BDD\ consists of the propositional connectives, the universal and existential quantifiers $\forall$ and $\exists$, a countable set of individual variables $\Var=\{x_1,x_2,\ldots\}$, a countable set of individual constants $\Const=\{c_1,c_2,\ldots\}$, and countable set of predicate symbols (each of a given finite arity) $\Pred=\{P_1,P_2,\ldots\}$; in this paper, we omit function symbols for the sake of simplicity. 
The set of formulas is defined as follows (in the Backus--Naur form):
\begin{equation*}
    \f\mathrel{::=}P(t_1,\ldots,t_n)\mid\bdnega\f\mid(\f\minc\f)\mid(\f\maxd\f)\mid\bddelta\f\mid(\forall x)\f\mid(\exists x)\f,
\end{equation*}
where $t_1,\ldots,t_n$ are terms, i.e., individual constants or variables. We employ the usual conventions for omitting parentheses and using the derived connectives.
The notions of free and bound variable and sentence are defined as usual.

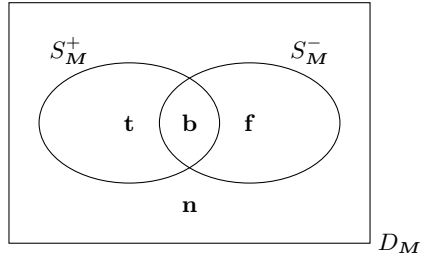
\begin{figure}\label{fig:ExtensionsExample}
    \centering
\begin{tikzpicture}[scale=0.8]
    \draw (0,0) rectangle (6,4);
    
    \draw (2,2) ellipse (1.5 and 1);
    
    \draw (4,2) ellipse (1.5 and 1);
    
    \node at (6.5,0) {$D_{\model M}$};
    \node at (1,3.2) {$S^+_{\model M}$};
    \node at (5,3.2) {$S^-_{\model M}$};
    \node at (2,2) {$\tct$};
    \node at (3,2) {$\tcb$};
    \node at (3,0.6) {$\tcn$};
    \node at (4,2) {$\tcf$};
\end{tikzpicture}
    \caption{The positive and negative extensions of a unary predicate $S$ in a model $\model M$.}
\end{figure}

A model for the (first-order) logic \BDD\ is a tuple
\begin{equation*}
    \model M=\bigl\langle
        D_{\model M},
        (P_{\model M})_{P\in\Pred},
        (c_{\model M})_{c\in\Const}
    \bigr\rangle,
\end{equation*}
where $D_{\model M}$ is a non-empty set (the \emph{domain} of the model) and $c_{\model M}\in D_{\model M}$ for all $c\in\Const$.
Each predicate $P_{\model M}\in\Pred$ of arity~$n$ is interpreted in~\model M by a function
\begin{equation*}
    P_{\model M}\colon (D_{\model M})^n\to \fclass{0,1}^2.
\end{equation*}
It can be decomposed into a pair of functions,
\begin{equation*}
    P^+_{\model M},P^-_{\model M}\colon (D_{\model M})^n\to \fclass{0,1},    
\end{equation*}
called the \emph{positive} and \emph{negative extension} (or the \emph{extension} and \emph{anti-extension}) of~$P$ in~\model M, so that for all $\vectn a\in D_{\model M}$,
\begin{equation*}
P_{\model M}(\vectn a)=\tuple{P^+_{\model M}(\vectn a),P^-_{\model M}(\vectn a)}.
\end{equation*}
For instance, if $S_{\model M}(a)=\tcb=\tuple{1,1}$, then $S^+_{\model M}(a)=S^-_{\model M}(a)=1$, for a unary predicate $S\in\Pred$ and $a\in D_{\model M}$ (see Figure~\ref{fig:ExtensionsExample}).

An \emph{evaluation} of individual variables in~\model M is a function
$e\colon\Var\to D_{\model M}$.
By $e[x\mapsto a]$, where $x\in\Var$ and $a\in D_{\model M}$,
we denote the evaluation $e'$ such that $e'(x)=a$ and $e'(y)=e(y)$ for each $y\in\Var$ different from~$x$.

The semantic value $\SemVal{t}Me\in D_{\model M}$ of a term~$t$ in a model~\model M under an evaluation~$e$ is defined as follows:
\begin{alignat*}{3}
    \SemVal{x}Me&=e(x)
        &\quad&\text{for each }x\in\Var
\\
    \SemVal{c}Me&=c_{\model M}
        &\quad&\text{for each }c\in\Const
\end{alignat*}
The truth value
$\SemVal{\f}Me=\bigl\langle\SemValPlus{\f}Me,\SemValMinus{\f}Me\bigr\rangle\in\fclass{0,1}^2$
of a formula~\f\ in a model~\model M under an evaluation~$e$ is given by the following Tarski conditions:
\begin{align}
\label{eq:BDD-predicate-Tarski-first}
    \SemVal{P(\vectn t)}Me &=
    P_{\model M}(\SemVal{t_1}Me,\dots,\SemVal{t_n}Me\bigr)
\\[1ex]
\SemVal{\heartsuit(\f_1,\dots,\f_m)}Me &=
    F_{\heartsuit}(\SemVal{\f_1}Me,\dots,\SemVal{\f_m}Me\bigr)
\\[1ex]
    \SemValPlus{(\forall x)\f}{M}{e} &= \inf_{a\in D_{\model M}}\SemValPlus{\f}{M}{e[x\mapsto a]}
\\
    \SemValMinus{(\forall x)\f}{M}{e} &= \sup_{a\in D_{\model M}}\SemValMinus{\f}{M}{e[x\mapsto a]}
\\[1ex]
    \SemValPlus{(\exists x)\f}{M}{e} &= \sup_{a\in D_{\model M}}\SemValPlus{\f}{M}{e[x\mapsto a]}
\\
\label{eq:BDD-predicate-Tarski-last}
    \SemValMinus{(\exists x)\f}{M}{e} &= \inf_{a\in D_{\model M}}\SemValMinus{\f}{M}{e[x\mapsto a]},
\end{align}
for each $n$-ary predicate $P\in\Pred$ and each $m$-ary connective~$\heartsuit$ of~\BDD\ (where the function $F_\heartsuit$ is given by the truth tables above).

The notion of logical consequence is defined as usual, i.e., as preservation of designated truth values from the premises to the conclusion. In more detail, for a given class \classmod K of models, we say that a set of \BDD-formulas $\mathrm\Gamma$ (positively) \emph{entails} a \BDD-formula~\f\ in~\classmod K, written $\mathrm\Gamma\models_{\classmod K}\f$, if for all models $\model M\in\classmod K$, the following condition holds: if $\SemValPlus{\p}Me=1$ for all $\p\in\mathrm\Gamma$ and all evaluations~$e$ in~\model M, then also $\SemValPlus{\f}Me=1$ for all evaluations~$e$ in~\model M.
We say that~\classmod K is the class of models of a \emph{theory}~$\mathrm\Gamma$ (i.e., a set of formulas called the \emph{axioms} of~$\mathrm\Gamma$) if $\SemValPlus{\f}Me=1$ for all $\model M\in\classmod K$, all evaluations~$e$ in~$\model M$, and all~$\f\in\mathrm\Gamma$.

By convention, we can drop the subscript in $\models_{\classmod K}$ if \classmod K is the class of all \BDD-models.
We may also omit the subscripts in $\SemVal\f Me$ if they are clear from the context.

\section{Four-Valued Free Quantification over the Logic \BDD}
\label{s:free-BDD}

As discussed in Section~\ref{s:intro}, we define free quantification over \BDD\ in a similar manner as has been done over the three-valued strong Kleene logic~\logic{K_3} (with truth-value gaps, cf.~\cite{Behounek-Dvorak:NondenotingTermsFuzzy-EUSFLAT}) and the logic \logic{LFI1} (with truth-value gluts,~\cite{Carnieli_FreeLFL1}), i.e., by introducing a primitive predicate of existence and restricting quantifiers to it.

The language of the (positive dual-domain) free logic over \BDD\ thus extends the language of \BDD\ just by the unary \emph{predicate of existence,} traditionally denoted by~\exP. Intuitively, \exP\ is intended to delimit the objects that exist.
For simplicity, in this paper we assume that \exP\ is bivalent, i.e., that all objects either do or do not exist, excluding the possibility of objects that both have and do not have existence, as well as objects that neither exist nor don't exist. This is ensured by assuming the following axiom for~\exP:
\begin{equation}
\label{eq:cons-E!}
    \bdcons\exP x
\end{equation}
This setting leads to a dual-domain semantics that considers two bivalent domains
in a given model~$\model M$:
\begin{itemize}
\item
    The \emph{outer domain} $\outD$ is the non-empty domain of the model, $(\outD)_{\model M}=D_{\model M}$.
\item
    The \emph{inner domain} $D_1$ is the positive extension of $\exP$, i.e., $(\inD)_{\model M}=\exP^+_{\model M}$.
\end{itemize}
In the following, we fix a model~\model M and can omit the subscripts specifying it.
In the inner domain~$\inD\subseteq\outD$ of the model, existent objects are collected. The elements of $\exP^-= \outD\setminus\inD$, i.e., the anti-extension of~$\exP$ in~\model M, serve as the absent referents of non-denoting terms, and thus represent nonexistent objects. 
The existence predicate makes it possible to restrict quantification to existing objects only, which leads to two kinds of quantifiers:
\begin{itemize}
\item
    The so-called \emph{outer quantifiers,} which range over the outer domain~\outD. These are the standard quantifiers of~\BDD, as they range over the whole domain $\outD=D$ of the model. In the dual-domain free logic, they are denoted by $\outAll,\outExi$ instead of $\inAll,\inExi$, as the latter symbols are reserved for the (more commonly used) inner quantifiers over the existing objects.

\item
    The \emph{inner quantifiers} $\inAll,\inExi$, which range over the inner domain~\inD, i.e., over the existing objects only.
\end{itemize}

The inner quantifiers can be defined from the standard (i.e., outer) \BDD-quantifiers $\outAll,\outExi$ by restricting them to the inner domain~\inD\ of existing objects, via relativization by the (bivalent) existence predicate~\exP:
\begin{eqnarray}
\label{eq:all1-crisp}
	(\inAll x)\f &\equivdef&
		(\outAll x)(\exP x \impl \f )
\\\label{eq:exi1-crisp}
	(\inExi x)\f &\equivdef&
    	(\outExi x)(\exP x \minc \f )
\end{eqnarray}
It can be observed that the semantics of the \BDD-connectives (given by the truth tables in Section~\ref{s:BDD}) ensures the intended semantics of the inner quantifiers, i.e., that they range over the inner domain $D_1=\exP^+$ only.
Indeed, the resulting Tarski conditions for the inner quantifiers read as follows:
\begin{align*}
\SemValPlus{(\forall x)\f}{M}{e} &= \inf_{a\in(D_1)_{\model M}}\SemValPlus{\f}{M}{e[x\mapsto a]}
\\
\SemValMinus{(\forall x)\f}{M}{e} &= \sup_{a\in(D_1)_{\model M}}\SemValMinus{\f}{M}{e[x\mapsto a]}
\\[1ex]
\SemValPlus{(\exists x)\f}{M}{e} &= \sup_{a\in(D_1)_{\model M}}\SemValPlus{\f}{M}{e[x\mapsto a]}
\\
\SemValMinus{(\exists x)\f}{M}{e} &= \inf_{a\in(D_1)_{\model M}}\SemValMinus{\f}{M}{e[x\mapsto a]}
\end{align*}

\begin{remark}
In practice, the most commonly used quantifiers are the inner ones, as the intended range of quantification is usually just over existing individuals: for instance, the proposition ``some horses fly'' is normally considered false, disregarding Pegasus and other fictitious flying horses; thus its adequate formalization is $(\inExi x)(Hx\minc Fx)$.
The outer quantifiers are needed for formalization of such propositions as ``some things do not exist'' (formalized as $(\outExi x)\bdnega\exP x$), which take into account nonexistent objects as well.
\end{remark}

\begin{remark}
Besides axiom~\eqref{eq:cons-E!} that ensures the bivalence of the existence predicate, further optional axioms can be adopted in free~\BDD. For example, it is reasonable to assume that although non-existing objects such as the round square can be inconsistent, the existing objects in the inner domain $D_1$ can only have non-contradictory (or even normal, i.e., \fclass{\tct,\tcf}-valued) properties. The normality on the inner domain can be ensured by the axioms
\begin{equation}\label{eq:D1-bivalent}
    (\inAll x_1)\dots(\inAll x_n)\,\bdcons P(\vectn x)
\end{equation}
for all predicates~$P$ in the language and the non-contradictoriness by the analogous schema
\begin{equation}\label{eq:D1-noncontradictory}
    (\inAll x_1)\dots(\inAll x_n)\,\bdnegad\bigl(P(\vectn x)\minc\bdnega P(\vectn x)\bigr)
\end{equation}
for all $P\in\Pred$ (where $n$ is the arity of~$P$).
\end{remark}

\begin{remark}
The described dual-domain free logic over~\BDD\ can be understood as a formal theory over~\BDD\ with the special predicate~\exP\ and special axiom~\eqref{eq:cons-E!} (stronger variants of the theory can also use axiom~\eqref{eq:D1-bivalent} or~\eqref{eq:D1-noncontradictory}), with the notational convention by which the quantifiers are written as~$\outAll$ and~$\outExi$, whereas the inner quantfiers and the derived connectives are regarded as abbreviations that can be expanded according to their definitions~\eqref{eq:all1-crisp}--\eqref{eq:exi1-crisp} and~\eqref{eq:def-bdnegad-BDD}--\eqref{eq:def-bdcons-BDD}. The known strong completeness theorem for~\BDD~\cite[Th.~4.9]{SanoOmori:BDDelta2014} then provides a strongly complete (natural deduction style) axiomatic system for the described variants of free logic over~\BDD.
\end{remark}

\begin{example}
As expected in free logic, the classical law of specification and its existential dual do not hold in free \BDD\ for the inner quantifiers:
\begin{align*}
    (\inAll x)\f(x)&\not\models\f(t),
&   \f(t)&\not\models(\inExi x)\f(x).
\\\intertext{However, they do hold if the existence of~$t$ is explicitly assumed:}
    \exP t,(\inAll x)\f(x)&\models\f(t),
&   \exP t,\f(t)&\models(\inExi x)\f(x).
\end{align*}

Moreover, while $\models\bigl(\outExi x)(\f(x)\maxd\bdnegad\f(x)\bigr)$, the same does not hold for the inner quantifier,
$\not\models(\inExi x)\bigl(\f(x)\maxd\bdnegad\f(x)\bigr)$; rather, the non-emptiness of the inner domain need be explicitly assumed:
$(\inExi x)\exP x\models(\inExi x)\bigl(\f(x)\maxd\bdnegad\f(x)\bigr)$.
\end{example}

\section{Free fuzzy logic with graded gaps and gluts}
\label{s:fuzzy}

We have intentionally formulated the definitions of the logic~\BDD\ and its dual-domain free variant in a format that can be easily fuzzified. In the fuzzification, the main change is replacing the bilattice $\fclass{0,1}^2$ of truth values with a bilattice~$L^2$, for a suitable residuated lattice~$L$; most definitions then stand as stated in the previous sections or require just minor adjustments. For simplicity, in this paper we will only consider the standard MV$\!_\baaz$-algebra $[0,1]_{\Luk_\baaz}$ of {\L}ukasiewicz fuzzy logic (with the operator~\baaz) as the underlying residuated lattice.
As the resulting logic uses the set $[0,1]^2$ of truth values, it is a bilattice-valued \emph{square fuzzy logic} (cf.~\cite{Turunen-Ozturk-Tsoukias:ParaconsistentPavelka,Genito-Gerla:BilatticeMultivaluedLogic}); we will denote the logic by \LBDD.

The set of truth values of the logic~\LBDD\ is $[0,1]^2$. The first component~\al\ of a truth value
$\tuple{\al,\be}\in[0,1]^2$ 
indicates the degree to which the proposition is true and the second component~\be\ indicates the degree to which it is false. Similarly to~\BDD, the degrees of truth and falsity of a proposition are evaluated independently, so unlike in standard fuzzy logics, they need not sum up to~$1$.
We keep the definitions of the truth constants $\tct=\tuple{1,0}$, $\tcf=\tuple{0,1}$, $\tcn=\tuple{0,0}$, and $\tcb=\tuple{1,1}$, which now refer to the corners of the square $[0,1]^2$.
The original truth values of fuzzy logic are embedded in $[0,1]^2$ as the pairs $\tuple{\al,1-\al}$ for $\al\in[0,1]$; we will call these values \emph{(fully) normal.} The truth values $\tuple{\al,\be}$ where $\al+\be<1$ are \emph{gappy} and those where $\al+\be>1$ \emph{glutty;} the former can be viewed as partial (or underdetermined) and the latter as contradictory (or overdetermined).
The \emph{designated} truth values of the logic~\LBDD\ are those true to the full degree, i.e., the pairs of the form $\tuple{1,\be}$ for all $\be\in[0,1]$. The information order~\ile\ and the truth order~\tle\ are defined on~$[0,1]^2$ as follows:
\begin{itemize}
\item
    $\tuple{\al_1,\be_1}\ile\tuple{\al_2,\be_2}$ iff $\al_1\le\al_2$ and $\be_1\le\be_2$;
\item
    $\tuple{\al_1,\be_1}\tle\tuple{\al_2,\be_2}$ iff $\al_1\le\al_2$ and $\be_1\ge\be_2$.
\end{itemize}

The notion of valuation is defined as in propositional~\BDD\ (see Section~\ref{s:BDD}).
The propositional language of~\LBDD\ contains the connectives of~\BDD\ with the addition of the connectives specific to the fuzzy logic~\LukD, namely, strong conjunction~\conj, strong disjunction~\disj, and the Delta operator~\baaz. The Tarski conditions for the primitive connectives of~\BDD\ are defined exactly as in~\eqref{eq:BDD-conn-Tarski-first}--\eqref{eq:BDD-conn-Tarski-last}, only evaluated in~$[0,1]$.
The Tarski conditions for the additional connectives are as follows:
\begin{align*}
  v^+(\f\conj\p)&=\max(v^+(\f)+v^+(\p)-1,0)\\   
  v^-(\f\conj\p)&=\min(v^-(\f)+v^-(\p),1)
\\[1ex]
  v^+(\f\disj\p)&=\min(v^+(\f)+v^+(\p),1)\\
  v^-(\f\disj\p)&=\max(v^-(\f)+v^-(\p)-1,0)
\\[1ex]
  v^+(\baaz\f)&=1-\sgn(1-v^+(\f))\\
  v^-(\baaz\f)&=\sgn(1-v^+(\f))
  \end{align*}

In the \LBDD-definitions of the derived connectives~\impl\ and~\bdcons, we must choose between the weak ($\minc,\maxd$) and strong ($\conj,\disj$) connectives of~\LBDD\ in place of the single conjunction and disjunction ($\minc,\maxd$) of~\BDD. In order for the definitions to conform with the intended semantics based on {\L}ukasiewicz logic, the strong connectives ($\conj,\disj$) are the appropriate choice. The \LBDD-definitions of the derived connectives thus read as follows:
\begin{align*}
    \bdnegad\f &\equiv
        \bdnega\bddelta\f
\\ 
    \f\impl\p &\equiv
        \bdnega\f\disj\p
\\ 
    \bdcons\f &\equiv
        (\bddelta\f\disj\bddelta\bdnega\f)\conj
        (\bdnega\bddelta\f\disj\bdnega\bddelta\bdnega\f)
\end{align*}
For the derived conncectives of~\LBDD, we obtain the following Tarski conditions:
\begin{align*}
  v^+(\bdnegad\f)&=1-v^+(\f)\\
  v^-(\bdnegad\f)&=v^+(\f)
\\[1ex]
  v^+(\f\impl\p)&=\min(v^-(\f)+v^+(\p), 1)\\
  v^-(\f\impl\p)&=\max(v^+(\f)+v^-(\p)-1, 0)
  \\[1ex]
  v^+(\bdcons\f) &= 
1-\bigl|v^+(\f)+v^-(\f)-1\bigr|\\
  v^-(\bdcons\f)&=1-v^+(\bdcons\f)
\end{align*}

\begin{figure}[ht]
\centering
\begin{tikzpicture}
  \draw[->] (-2.5,0) -- (4,0) node[right] {$D_{\model M}$};
  \draw[->] (-2,-0.5) -- (-2,1.2) node[left] {$1$};
  \node[left] at (-2, -0.2) {$0$};
    
  \draw[domain=-2:3.5,smooth,variable=\x,blue] plot ({\x},{exp(-\x*\x)});
  \draw[dashed] (-2,0) -- (-2,{exp(-9)});
  \draw[dashed] (3.5,0) -- (3,{exp(-9)});
  \node[above] at (-0.5, 1) {$P^+_{\model M}$};

  \draw[domain=-2:3.5,smooth,variable=\x,red] plot ({\x},{exp(-(\x-1)*(\x-1)/1.2)});
  \draw[dashed] (-2,0) -- (-2,{0.7*exp(-7.2)});
  \draw[dashed] (3.5,0) -- (3,{0.7*exp(-7.2)});
  \node[above] at (1.5, 1) {$P^-_{\model M}$};
\end{tikzpicture}
\caption{The positive and negative extensions of a unary predicate $P$ in a model $\model M$ of~\LBDD.}
\label{fig:extensions-fuzzy} 
\end{figure}
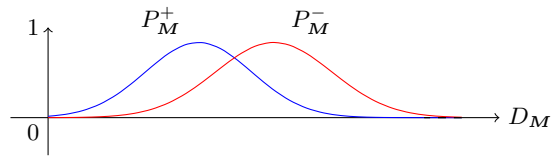

The first-order language of \LBDD\ is the same as that of \BDD, modulo the added connectives ($\conj,\disj,\baaz$) of~\LukD. 
The definition of a model for \LBDD\ is also the same as for~\BDD, with the only difference that predicate symbols are evaluated in $[0,1]^2$ instead of $\fclass{0,1}^2$, i.e.,
\begin{equation*}
        P_{\model M}\colon (D_{\model M})^n\to[0,1]^2.
\end{equation*}
Just like in \BDD, $P_{\model M}$ can be decomposed into its positive and negative extensions $P^+_{\model M},P^-_{\model M}\colon (D_{\model M})^n\to[0,1]$, which are ordinary $[0,1]$-valued fuzzy sets or $n$-ary fuzzy relations on~$D_{\model M}$ (see Figure~\ref{fig:extensions-fuzzy}), so that
\begin{equation*}
P_{\model M}(\vectn a)=\tuple{P^+_{\model M}(\vectn a),P^-_{\model M}(\vectn a)}\in[0,1]^2.
\end{equation*}

The notion of evaluation and the semantic value of a term in a model are defined in exactly the same way as in \BDD.
Also the truth value of a formula, $\SemVal{\f}Me=\bigl\langle\SemValPlus{\f}Me,\SemValMinus{\f}Me\bigr\rangle$, is defined by the same Tarski conditions~\eqref{eq:BDD-predicate-Tarski-first}--\eqref{eq:BDD-predicate-Tarski-last}, only now valued in $[0,1]^2$ instead of $\fclass{0,1}^2$.
Finally, the (positive) consequence relation and the class of models of a theory in \LBDD\ are also defined by the same conditions as given in Section~\ref{s:BDD} for~\BDD.

The notions of the dual-domain free logic over~\LBDD, i.e., the existence predicate~\exP, the outer domain $(D_0)=D_{\model M}$, the inner domain $(D_1)_{\model M}={\exP}^+_{\model M}$, and the outer quantifiers $\outAll,\outExi$ are also defined as in Section~\ref{s:free-BDD} over \BDD.
In \LBDD\, the role of axiom~\eqref{eq:cons-E!} is different, as it only ensures the full normality of the existence predicate but allows the inner domain $(D_1)_{\model M}={\exP}^+_{\model M}$ to be fuzzy; in order to ensure the bivalence of the inner domain, the axiom needs to be adjusted as follows:
\begin{equation}
\label{eq:exP-bivalent-LBDD}
    \bdcons\exP x\minc(\exP x\maxd\bdnega\exP x)
\end{equation}
Assuming~\eqref{eq:exP-bivalent-LBDD}, the definitions of the inner quantifiers can use the same formulas~\eqref{eq:all1-crisp}--\eqref{eq:exi1-crisp} as in~\BDD. Like in the dual-domain semantics over~\BDD, it is reasonable to assume the full normality of all predicates on the inner domain, which is ensured by the axiom~\eqref{eq:D1-bivalent}.

\begin{remark}
It can be observed that adding the following axiom schema to~\LBDD\ reduces it (including the dual-domain semantics) to the four values of~\BDD:
\begin{equation*}
    \bigl(\bdcons\f\minc(\f\maxd\bdnega\f)\bigr)\maxd\bdnega\bdcons\f
\end{equation*}
Similarly, the axiom schema
$\bdcons\f\maxd(\bdnega\bddelta\f\minc\bdnega\bddelta\bdnega\f)$
reduces the setting to that of partial fuzzy logic from~\cite{Behounek-Dvorak:NondenotingTermsFuzzy-EUSFLAT}, which considers only the value~\tcn\ (denoted there by~$*$)
besides the degrees from~$[0,1]$.
Using the axiom schema
$\bdcons\f\maxd\bdnega\bdcons\f$
in~\LBDD\ adds the value~\tcb\ to the setting of~\cite{Behounek-Dvorak:NondenotingTermsFuzzy-EUSFLAT}, giving rise to a logic with the truth values $[0,1]\cup\fclass{\tcn,\tcb}$.
\end{remark}

\section{Conclusion}
\label{s:conclusion}

In this paper, we have presented a four-valued system of dual-domain positive free logic with truth-value gaps and gluts over the logic~\BDD, sketched a fuzzification of~\BDD\ via {\L}ukasiewicz fuzzy logic~\LukD, and defined the dual-domain semantics over the resulting logic~\LBDD.
While the free logic over~\BDD\ is axiomatized by the strong completeness theorem for~\BDD\ due to Sano and Omori~\cite{SanoOmori:BDDelta2014}, we have only presented the (first-order and dual-domain) semantics of~\LBDD; a sound and complete axiomatization of~\LBDD\ is left for future work.

It can be observed that free logics based on truth-functional underlying logics are susceptible to lottery-style paradoxes, where disjuncts are evaluated by~\tcn, while their exhaustive disjunction should have the value~\tct, contrary to the semantics of \logic{K_3} or \logic{(\Luk)\BDD}. A remedy can take the form of using a non-truth-functional modality to handle such cases, similar to the two-layered modalities introduced in~\cite{Bilkova-et-al:QualitativeTwoLayered,Bilkova-et-al:TwoLayeredParaconsistentProbabilities} over related logics. This option is also left for future work. 

\paragraph*{Acknowledgment.}
L.~B\v{e}hounek and A.~Dvo\v{r}\'ak were supported by project No.~22-01137S of the Czech Science Foundation.

%
%
%


\end{document}